\documentclass[sn-mathphys]{sn-jnl}
\jyear{2024}%

\usepackage{upquote}
\theoremstyle{thmstyleone}%
\newtheorem{theorem}{Theorem}

\def\Im{\hbox{\rm Im\kern .7pt}}
\def\Re{\hbox{\rm Re\kern .7pt}}
\def\complex{{\mathbb C}}
\def\real{{\mathbb R}}
\def\qed{\hbox{\vrule width 2.5pt depth 2.5 pt height 3.5 pt}}
\def\too{\,\to\,}

\theoremstyle{thmstylethree}%
\newtheorem{definition}{Definition}%
\def\half{\textstyle{1\over 2}}

\begin{document}

\title[Numerical computation of the Schwarz function]{Numerical computation of the Schwarz function}

\author*{\fnm{Lloyd N.} \sur{Trefethen}}\email{trefethen@seas.harvard.edu}

\affil*{\orgdiv{School of Engineering and Applied Sciences},
\orgname{Harvard University},
\orgaddress{\street{29 Oxford St.}, \city{Cambridge, MA}, \postcode{02138},
\country{USA}}}

\abstract{An analytic function can be continued across
an analytic arc $\Gamma$ with the help of
the Schwarz function $S(z)$, the analytic function 
satisfying $S(z) = \bar z$ for $z\in \Gamma$.
We show how $S(z)$ can
be computed with the AAA algorithm of rational 
approximation, an operation that is the basis of 
the AAALS method for solution of Laplace and
related PDE problems in the plane.
We discuss the challenge of
computing $S(z)$ further away from from $\Gamma$, where
it becomes multi-valued.}



\maketitle

\section{Introduction}\label{intro}

The Schwarz function for an analytic arc $\Gamma$ is a beautiful
and long-established idea.  It is simply the analytic function
$S(z)$ defined near $\Gamma$ that takes the values $S(z)
= \bar z$ for $z\in \Gamma$, and it holds the key to the
continuation of analytic functions across $\Gamma$, as well as
the continuation of related functions such as harmonic potentials
and Helmholtz wavefields.  The ideas go back in part to Schwarz
in 1870~\cite[pp.~150--151]{schwarz}, Grave in 1895~\cite{grave}, and Herglotz
in 1914~\cite[p.~305]{herglotz}, and the function was given its name by
Davis and Pollak in 1958~\cite{dp}.  Two books have been published,
by Davis~\cite{davis} and Shapiro~\cite{shapiro}.  Yet the Schwarz
function is not a widely known or widely used tool, and a reason
for this may be that until recently, there has been no method
available for computing it numerically, so it has played no role
in numerical algorithms.  Here we introduce such a method, based
on AAA rational approximation.  Examples of the method were shown
earlier in~\cite{analcont} and~\cite{nonconvex}.

\section{The Schwarz function}\label{schwarzfun}

The starting point of analytic continuation is the Schwarz reflection
principle, whose simplest case is sketched on the left side of
Figure 1.  Let $\Omega\subseteq \complex$ be a connected open set
that is symmetric about $\real$, $\overline{\Omega} = \Omega$,
and whose intersection with $\real$ is a finite or infinite open
interval $\Gamma$.  (Here $\overline{\Omega}$ means the complex
conjugate of $\Omega$, not its closure.)  We define $\Omega^+ =
\{z\in \Omega: \Im (z) > 0\}$ and $\Omega^- = \{z\in \Omega: \Im
(z) < 0\}$.

\begin{theorem}
{\bf Schwarz reflection principle.} Let\/ $\Omega$, $\Gamma$,
$\Omega^+$ and\/ $\Omega^-$ be as defined above, and let $f$
be analytic in the interior of\/ $\Omega^+$ and continuous
on $\Omega^+\cup \Gamma$, taking real values on $\Gamma$.
Then the formula \begin{equation} f(\bar z ) = \overline{f(z)}
\label{analcont} \end{equation} defines an analytic continuation of\/
$f$ to all of\/ $\Omega$.  \end{theorem}

{\em Proof.}  Let $f$ denote the function defined in all of\/
$\Omega$ by (\ref{analcont}), which makes sense as the two sides
of the equation agree for $z\in \Gamma$ since $z$ and $f(z)$ are
real there.  It is immediate that $f$ is analytic in the interiors
of $\Omega^+$ and $\Omega^-$ (the Taylor series about two points $z$
and $\bar z$ are the same apart from conjugation of the coefficients)
and that $f$ is continuous throughout $\Omega$.  It remains only
to show that $f$ is analytic on $\Gamma$.  This follows from a general
principle involving analytic functions on adjacent domains that take
the same continuous limit on an interface curve, which can
be proved by a contour integral argument~\cite[p.~183]{nehari}. \qed

\begin{figure}
\begin{center}
\vspace{10pt}
\includegraphics[clip, scale=.75]{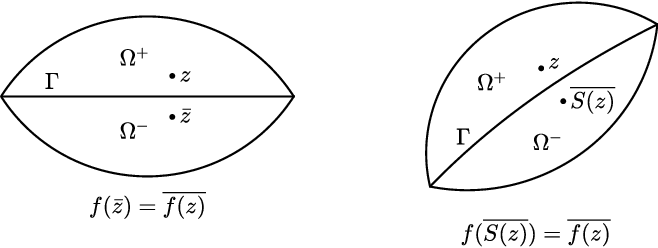}
\vspace{4pt}
\end{center}
\caption{\small On the left, the Schwarz reflection principle
(Theorem 1) extends an analytic function that takes real values
on a real interval across the interval from one half-plane to
the other.  On the right, the generalization by the Schwarz function
$S(z)$ for reflection of an analytic function
across an analytic arc $\Gamma$, again assuming it takes real
values on $\Gamma$.}
\end{figure}

\vskip 5pt

The Schwarz reflection principle is widely known and
has applications to many situations.
With the Schwarz function, we can generalize it to cases where
$\Gamma$ is not a real segment but an analytic curve,
as sketched on the right side of Figure 1.
Let $\Gamma$ be an analytic Jordan arc or Jordan curve,
by which we mean
the image in $\complex$ of the unit circle or an arc of the unit circle
under a one-to-one analytic function with nonvanishing derivative.
Consider the function
\begin{equation}
S(z) = \bar z,  \quad z\in \Gamma.
\label{Sdef}
\end{equation}
(This may seem like the archetype of a nonanalytic function,
but on an analytic curve as opposed to a domain with interior,
$z \mapsto\bar z$ is analytic.)  Like any other analytic function defined
on a curve in $\complex$, $S$ can be analytically continued
to a complex neighborhood of $\Gamma$, and this is the Schwarz
function.  In a word: $S(z)$ is the complex analytic function in a
neighborhood of $\Gamma$ that takes the values $S(z) = \bar z$ for
$z\in\Gamma$.  The point of $S$ is that $z\mapsto \overline{S(z)}$
is the generalization to a general arc $\Gamma$ of the reflection
$z\mapsto \bar z$ when $\Gamma$ is real.

We now explain how the analytic continuation works.  Since $S$ is analytic,
it must extend to an analytic function on some neighborhood
of $\Gamma$.  Since $|d\kern 1pt \overline{S(z)}/dz| = 1$ on $\Gamma$, it follows
that in a sufficiently small neighborhood of $\Gamma$, $\overline{S(z)}$ maps a point
$z$ close to $\Gamma$ on one side to a point close to $\Gamma$ on
the other side.  This function
$\overline{S(z)}$ is co-analytic, not analytic, but if we compose it with
itself we get an analytic function---the identity:
\begin{equation}
\overline{S( \overline{S(z)} )} = z.
\label{ident}
\end{equation}
To verify that (\ref{ident}) holds for $z$ in a neighborhood of
$\Gamma$, we note first that (\ref{Sdef}) implies that it holds
for $z\in\Gamma$.  It follows that on any neighborhood $B$ of $\Gamma$
such that $S$ can be continued to $B$ and then to $\overline{S(B)}$,
the equality persists.  Thus we see that the map $z\mapsto
\overline{S(z)}$ is indeed a reflection sufficiently
near $\Gamma$ in the sense that its square
is the identity.  Such a map is called an involution.  We make the
following definition, sketched on the right side of Figure 1.

\vskip 5pt
\begin{definition}
Given an analytic Jordan arc or Jordan curve
$\Gamma$, a {\em reflection domain} for $\Gamma$
is an open set $\Omega\subseteq\complex$ containing $\Gamma$
that is divided into two disjoint open sets $\Omega^-$
and $\Omega^+$ by $\Gamma$ and in which a single-valued
analytic Schwarz function $S$ for $\Gamma$ can
be defined that satisfies
\begin{equation}
\overline{S(\Omega^+)} = \Omega^-, \quad
\overline{S(\Omega^-)} = \Omega^+,
\label{defn1}
\end{equation}
and
\begin{equation}
\overline{S( \overline{S(z)} )} = z, \quad z\in \Omega.
\label{defn2}
\end{equation}
\end{definition}
\vskip 5pt

We can now state the generalization of the reflection principle.

\begin{theorem}
{\bf Reflection across an analytic curve.}
Let\/ $\Gamma$ be an analytic Jordan arc or Jordan curve,
let $\Omega$ be a reflection domain for $\Gamma$,
and let $f$ be analytic in the interior
of\/ $\Omega^+$ and continuous on
$\Omega^+\cup \Gamma$, taking real values on\/ $\Gamma$.
Then the formula
\begin{equation}
f(\overline{S(z)}) = \overline{f(z)}
\label{Scont}
\end{equation}
defines an analytic continuation of\/ $f$ to all of\/ $\Omega$.
\end{theorem}

{\em Proof.}  
The same as for Theorem 1, with the obvious changes. \qed

\vskip 5pt

Note that in both Theorems 1 and 2, we have supposed that
$f$ takes real values on $\Gamma$.  If it takes
imaginary values on $\Gamma$, the reflection formula adjusts to 
\begin{equation}
f(\overline{S(z)}) = \overline{-f(z)},
\label{Scontimag}
\end{equation}
and there are similar elementary adjustments if $\Gamma$
is another straight segment or a circular arc in $\complex$.
More generally, the situation in which $f$ takes complex values along an
analytic arc can be treated with the help of a second Schwarz
function~\cite[Prop.~1.3]{shapiro}, but we shall not consider
this case.

\vskip 5pt

{\em Example: unit circle.}  If $\Gamma$ is the unit circle, then
$S(z) = z^{-1}$.  Here we can take any annulus $r < |z| < r^{-1}$
with $r\in [\kern .5pt 0,1)$ as a reflection domain.  The choice
$r=0$ provides a reflection domain that is as large as possible,
the punctured plane $0 < |z|<\infty$, but for choices
of $\Gamma$ other than a line or a circle,
there may not be a maximal choice of $\Omega$, because
$S$ will usually have branch points, as in the next example.

\vskip 5pt

{\em Example: ellipse.}
For any $\rho>1$, define the {\em
$\rho$-ellipse} $E_\rho$ to be the image in the $z$-plane of the
circle $|w|=\rho$ in the $w$-plane under the Joukowsky map $z =
(w+w^{-1})/2$.  Geometrically, $E_\rho$ is the ellipse with foci
$\pm 1$ whose semiminor and semimajor axis lengths sum to $\rho$.
The Schwarz function for $E_\rho$ is~\cite[p.~25]{davis}
\begin{equation}
S(z) = \half(\kern 1pt\rho^2+\rho^{-2})
\kern .5pt z - \half (\kern 1pt\rho^2-\rho^{-2})\sqrt{z^2-1}.
\label{SE}
\end{equation}
Note that $S$ has a pair of branch points at $\pm 1$, and that
otherwise, it is free of singularities in $\complex$.  As a
reflection domain $\Omega$ for $S$, the obvious choice is to
take $\complex$ minus a branch cut $[-1,1]$, and also removing the
reflection of $[-1,1]$ in the ellipse, namely the segments $(-\infty,-a]$ and $[a,\infty)$
with $a = \half(\kern 1pt\rho^2+\rho^{-2})$.  However, other choices
of branch cuts and reflection domains are also valid.

\vskip 5pt

Besides the circle and ellipse, certain other domains have known Schwarz functions, such
as parabolas and hyperbolas, as reported in~\cite{davis}.  However, there are
not many such cases.  A good deal of attention has been given to the study of
so-called {\em quadrature domains}, which are domains for which the interior
has a Schwarz function that is meromorphic, i.e., analytic apart from
poles~\cite{gust,shapiro}.
These too are very special cases, as the nearly
universal appearance of approximate branch cuts in caclulations
like those illustrated in Figures 2--5 of this paper testifies.

\section{AAA approximation of the Schwarz function}\label{AAA}
In 2018 an algorithm for numerical rational approximation was published
that is a step change from previously available algorithms, making it now possible
to compute near-optimal approximations of
many functions on many domains in a fraction of a second on a
laptop.  This is the AAA algorithm, whose name derives from
``adaptive Antoulas-Anderson.''  The algorithm is a greedy iteration
based on a barycentric representation of a rational function $r$ with
coefficients chosen at each iterative step via solution of a linear
least-squares problem.  We will not give details, which can be found
in~\cite{aaa}.  AAA has implementations in Matlab~\cite{chebfun},
Python~\cite{virtanen},
and Julia~\cite{driscoll,macmillen}, and its
invocation in Matlab starts from the command
\begin{equation}
\hbox{\tt r = aaa(F,Z)}
\label{command}
\end{equation}
where $Z$ is a vector of sample points, typically a few hundreds or thousands,
and $F$ is a vector of corresponding function values.  The output is a
Matlab function handle corresponding to a rational function $r$ satisfying
\begin{equation}
{\| f - r\|\over \|f\|} \le 10^{-13}
\label{tolerance}
\end{equation}
in the default mode of operation.  Here $\|\cdot\|$ is the
$\infty$-norm over $Z$, which is typically an approximation to the
$\infty$-norm over a real or complex continuum of interest.  The AAA
rational approximation is not quite optimal for its degree, but its
error is usually within a factor on the order of $10$ of optimal.
Another way to say it is that the AAA approximant often has degree
about $10\%$ greater than that of the minimal-degree approximation
$r$ satisfying (\ref{tolerance}).

\begin{figure}[t!]
\begin{center}
\vspace{10pt}
\includegraphics[clip, scale=.80]{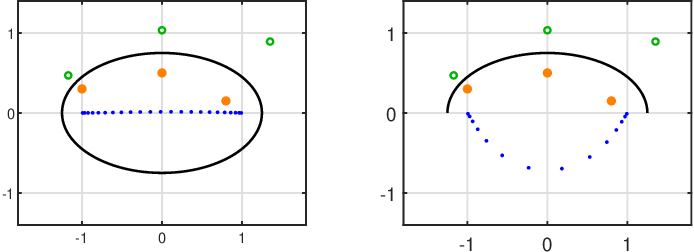}
\vskip +10pt
\includegraphics[trim=0 160 0 20, clip, scale=.80]{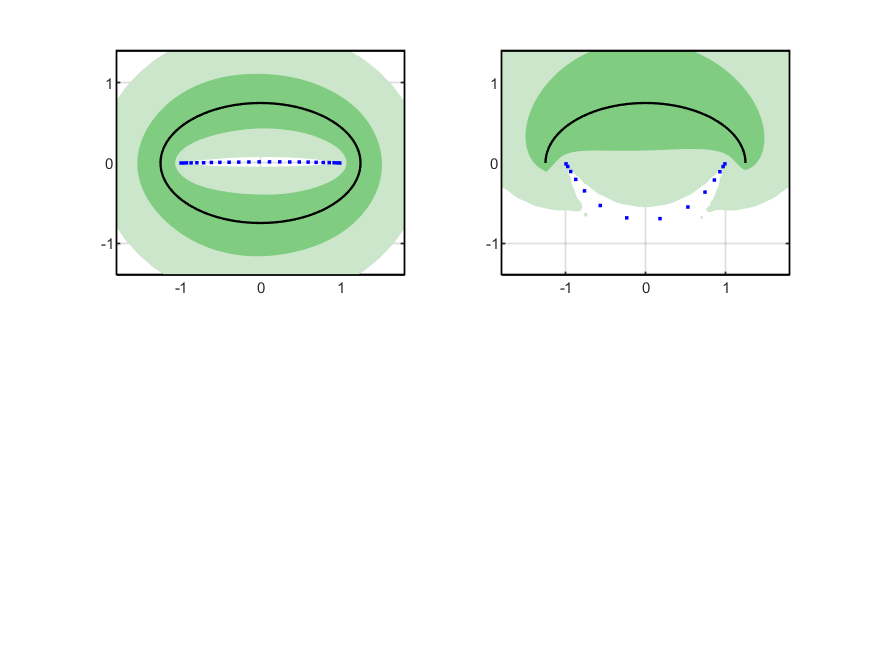}
\vspace{-14pt}
\end{center}
\caption{\label{figellipse}\small In the first row, poles (blue)
of AAA approximations to $S(z)$ for an ellipse and
a half-ellipse.
These computations took less than $0.01$ s on a laptop.
Both cases show a string of poles
approximating a branch cut extending from $-1$ to $1$.
In each image three
orange points $z$ have been selected and their
reflections $\overline{r(z)}\approx\overline{S(z)}$ shown as
green circles.
In the second row, error indicators $(\ref{measure})$ of the approximations
$r(z)\approx S(z)$.  In the dark green region the accuracy is better than
$10^{-8}$, and in the light green region it is better than $10^{-1}$.}
\end{figure}

Our proposed method for computing $S(z)$ is simply to apply the AAA
algorithm to a discrete set $Z$ of sample points along $\Gamma$,
with function values $F = \overline{Z}$.  The result is a rational
function $r$ that matches $S$ to $13$ digits of relative accuracy
on $\Gamma$.  (With extended-precision AAA implementations, such
as are available in Julia, one can get more digits of accuracy.)
The question then becomes, as $z$ moves away from $\Gamma$, how
much accuracy in $r(z)\approx S(z)$ is retained and for how far?

What almost always emerges from such computations is that $r$ has
strings of poles that approximate branch cuts, a phenomenon familiar
in rational approximation and analyzed for Pad\'e approximation by
Stahl~\cite{stahl97,stahl12} and for multipoint Pad\'e approximation
by Buslaev~\cite{bus1}.  As our first illustration, the top-left panel
of Figure 2 shows the poles of the AAA approximation to $S(z)$
for the case in which $\Gamma$ is the $\rho$-ellipse with $\rho=2$,
discretized by 100 points.  We see 23 poles lining up in $[-1,1]$,
nicely approximating a branch cut there.  The figure also shows
three arbitrarily chosen points $z$ near $\Gamma$ as orange dots and
their reflections $\overline{r(z)} \approx \overline{S(z)}$ as
green circles.  The top-right panel of the figure repeats the same computation,
except with $\Gamma$ taken as just the upper half of the ellipse,
which in principle has exactly the same Schwarz function.  The same
branch points are in evidence, but a different choice of branch cut.
In addition to the 14 poles visible in the plot, there is also a
pole off-scale at $\approx (4.4-7.4i)\times 10^{7}$, a numerical
approximation to $\infty$.

Where and how accurately does $r$ approximate $S$?    As one indicator,
following (\ref{defn2}), we can compute $\overline{r(\overline{r(z)})}$ and
see how close it comes to $z$ itself at various points $z\in\complex$.
The bottom row of Figure 2 applies this idea by filling in regions bounded by contours
\begin{equation}
|\overline{r(\overline{r(z)})} - z | = C
\label{measure}
\end{equation}
with $C = 10^{-8}$ (dark green) and $C = 10^{-1}$ (light
green).\footnote{I am
grateful to Keaton Burns for suggesting this method of visualization.}
We will have more to say about these images in section~5.

\begin{figure}[t!]
\begin{center}
\vspace{10pt}
\includegraphics[clip, scale=.85]{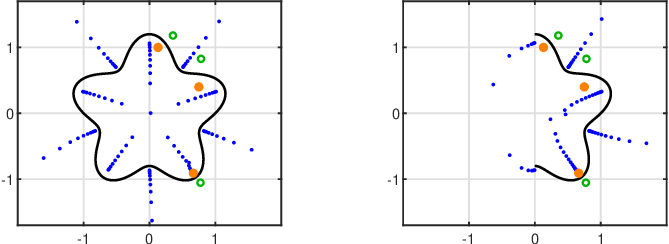}
\vskip 10pt
\includegraphics[clip, scale=.85]{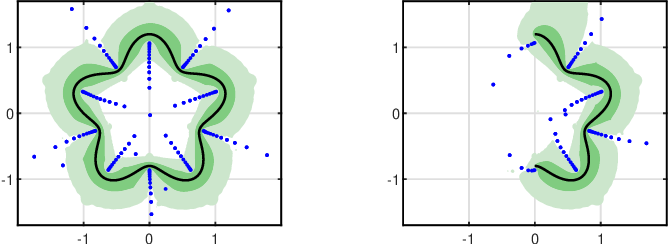}
\vspace{4pt}
\end{center}
\caption{\label{figsquiggle}\small Repetition of Figure~$\ref{figellipse}$ for
a different pair of curves $\Gamma$.}
\end{figure}

Figure 3 applies the same computation to a pair of curves $\Gamma$ for which
$S(z)$ is not known analytically, defined in polar coordinates by a radius
equal to $1 + 0.2\sin(5\theta)$.
Five branch cuts inside $\Omega$ and five more outside are in
evidence, with a convincing reflection domain between them.

\section{Singular and near-singular contours}
Figure 4 shows four more examples of computed approximations to
Schwarz functions.  In the first row, $\Gamma$ is an 
analytic curve as in our previous examples, whereas in
the second row it has $1$ and $6$ singular points in the
left and right images, respectively.

\begin{figure}
\begin{center}
\vspace{5pt}
\includegraphics[clip, scale=.85]{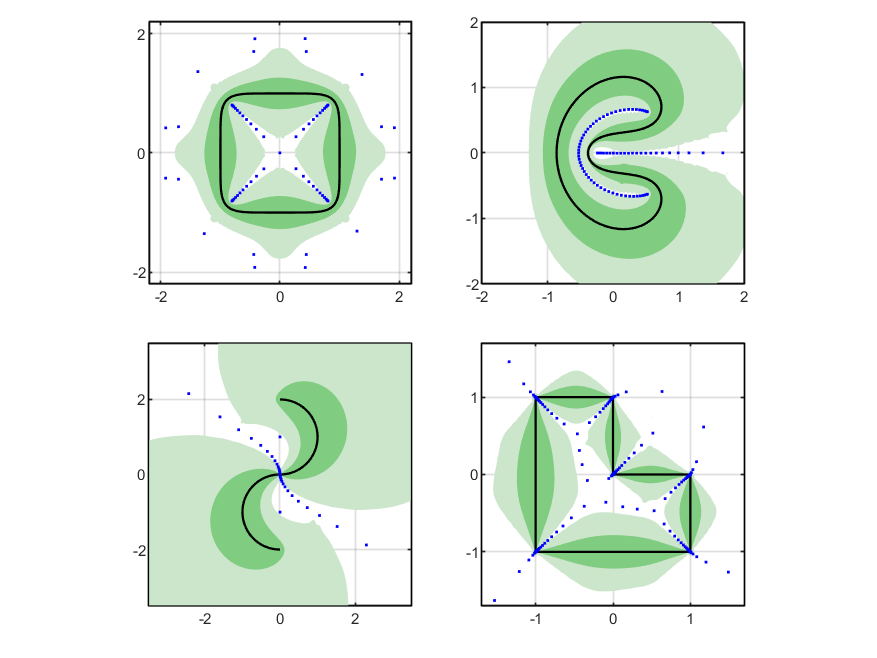}
\vspace{-12pt}
\end{center}
\caption{\label{figfour}\small Four more examples of numerically computed Schwarz functions.
In the third image $\Gamma$ has a singularity at $z=0$, and in the fourth one, there are six
singularities at the six corners.  In both of these cases there are several mathematically
independent branches of $S(z)$, which the rational approximations separate by
strings of poles approximating branch cuts.  For the final image the tolerance has been loosened
to $10^{-8}$ and the dark green region shows accuracy $10^{-5}$ rather than the
usual $10^{-8}$.}
\end{figure}

In the first image, $\Gamma$ is a superellipse defined by the
condition $x^6 + y^6 = 1$, discretized by 200 points.
(In~\cite{davis} and~\cite{dp} an explicit formula for $S(z)$
is given for the curve defined by $x^4 + y^4 = 1$.)
Besides the pole at the origin there are 68 further poles with
four-fold symmetry (four of them off-scale in this image).  The poles
line up along a branch cut in the form of a cross with
branch points approximately at $\pm 0.81 \pm 0.81i$.

The second image shows a domain with an inlet discretized
in 300 points, a configuration
where rational approximations are exponentially 
more powerful than polynomial approximations~\cite{nonconvex}
(exponentially with respect to the inlet length-to-width ratio).
There are 82 poles, six of them off-scale at
$z \approx 2.1,$ $2.9,$ $4.3,$ $7.5,$ $17.5,$ and $89.5$.

The third image introduces a singularity.  Here $\Gamma$ consists of an upper
semicircle joined to a lower semicircle, so that at the junction at $z=0$,
$\Gamma$ is $C^1$ but not $C^2$.  Each semicircle is discretized
by 400 points exponentially clustered near $z=0$, since
poles will need to cluster there.  Mathematically, there is one Schwarz
function for the upper half and another for the lower half, both of them
trivial since these are just semicircles.    If we wish to think of a single
global Schwarz function, then it will have two branches that are entirely independent,
with no branch point linking the two.  A global AAA rational
approximation, however, is forced
to choose a single-valued approximation, and one can see a curve of
68 poles dividing the semicircles (four of them off-scale)
in addition to the two poles very close to $\pm \kern .5pt i$.
For this computation, which essentially involves a complex sign function,
we found that for good results it was important
to run the AAA algorithm with the \verb|'sign'| modification described
in~\cite{zolo}.

The fourth image shows an L-shaped region, each of whose six sides
has been discretized by 300 points exponentially clustered at both ends.
Whereas the other calculations each take on the order of $0.1$ s on
a laptop, this one requires 7 seconds because of the large number of poles,
361.  In this configuration $\Gamma$ is $C^{\kern .6pt 0}$ but not $C^1$, putting it in
the class of problems investigated by Newman in his study of rational approximation
of $|x|$ on $[-1,1]$~\cite{newman}.  For such problems the convergence rate is
root-exponential, i.e., with errors decreasing at the rate $\exp(-C\sqrt n
\kern 1pt)$ as a function of degree $n$ for some $C>0$, and in this example
we loosened the AAA tolerance from $10^{-13}$ to $10^{-8}$.  The contour
level for the dark green region in the plot was correspondingly loosened
from the usual $10^{-8}$ to $10^{-5}$.  Mathematically, the Schwarz function
could be regarded as having six independent branches with no branch points, and
again the global rational approximation is forced to introduce approximate
branch cuts, as is evident in the figure.

\section{Multiple branches}
In Figures 2--4, both the dark and light green regions can be interpreted as
reflection domains as defined in section 2: to a certain numerical accuracy, the
approximations $\overline{r(z)}\approx \overline{S(z)}$ reflect these domains across $\Gamma$
into themselves.  

The question presents itself, what can be said about the white regions in these plots?
The first temptation is to suppose that these mark regions where the numerical approximation
has simply failed to capture $S(z)$.  However, the 
true situation is often better than this.
Because of branch points, $S(z)$ is a multivalued function,
and $r(z)$ is sometimes taking us
with quite good accuracy to another branch.

We can see this effect in Figure 2, the example of the ellipse, where everything is
known analytically.  Suppose we start at the point
$z=1.3\kern .6pt i$, which lies in the light green region,
and apply $\overline{r(z)}$ several times.
These are the values that result, showing just
the alternation we hope for as $z$ is reflected
back and forth across $\Gamma$:
\begin{equation}
1.3\kern .6pt i \too 0.3127i \too 1.3000\kern .6pt i \too 0.3127i \too 1.3000\kern .6pt i \too \cdots .
\end{equation}
On the other hand if we start at $z=3\kern .6pt i$, which lies
in the white region, this is what happens:
\begin{equation}
3\kern .6pt i \too -0.4457i \too -1.1057i \too -0.4457i \too -1.1057i \too \cdots .
\end{equation}
The initial reflection to $-0.4457i$ is extremely accurate, matching the
exact value from (\ref{SE}) to 14 digits.  Since the imaginary part is
negative, however, this value has crossed the approximate branch cut of
$r$.  When $r$ is applied to this new value, the number $-1.1057i$ matches
the correct value for the {\em other\/} branch of $S$ to 9 digits, as one
can confirm by applying (\ref{SE}) with the sign of the square root negated.
As a measure of this accuracy, Figure 5 repeats the lower half of Figure 2,
but now, instead of measuring the error by (\ref{measure}),
the plot colors regions according to
\begin{equation}
\min\{ \,|r(z) - S_1(z)|\kern 1pt,\kern 1pt |r(z) - S_2(z) |\,\} = C,
\label{measure2}
\end{equation}
where $S_1$ and $S_2$ are the two branches of the
Schwarz function as given in (\ref{SE}).
Now the blue ink covers a large area, except near the poles.
\begin{figure}
\begin{center}
\vspace{10pt}
\includegraphics[trim=0 160 0 20 clip, scale=.80]{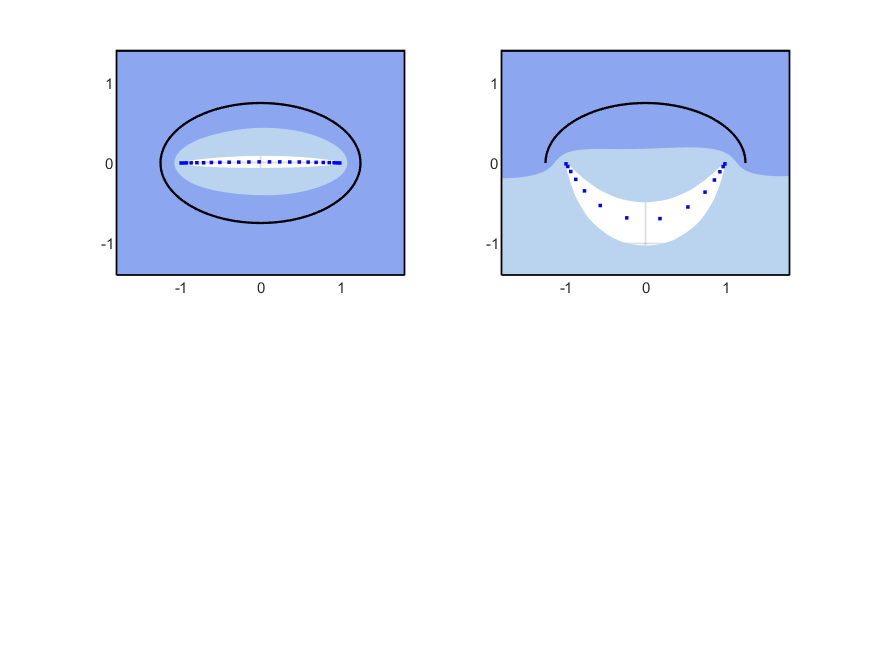}
\vspace{-14pt}
\end{center}
\caption{\label{figellipse3}\small Repetition of
the bottom row of Figure~\ref{figellipse} with
(\ref{measure}) replaced by (\ref{measure2}), so that
the plot indicates deviation at each point from the closer of
the two branches of the exact Schwarz function $S(z)$
of $(\ref{SE})$, and with colors
changed from green to blue.  Except along the
approximate branch cut, $r$ matches one or the other branch of $S(z)$
with good accuracy everywhere.  The dark and light regions extend distances
from the origin about $4000$ and $1.5\times 10^7$ on the left,
$6$ and $5000$ on the right.}
\end{figure}

We thus see that AAA rational approximations $r(z)$ contain information
that is easily interpretable in the green regions of Figures 2--4, and may also
contain accurate information in parts of the white regions.  Whether this
information can be usefully utilized, or perhaps improved upon by the
application of other computational methods, remains to be seen.
There is a small literature on numerical methods for tracking multivalued
analytic functions, but it is not very advanced~\cite{fhsw,loglight,shafer}.

\section{Discussion}
My first motivation for computing the Schwarz function concerns
the AAA-least squares (AAALS) method for solving 2D Laplace problems,
which was invented by Stefano Costa~\cite{costa}.
For the basic example, consider the Dirichlet problem
\begin{equation}
\Delta u(z) = 0, ~ z\in \Omega, \quad 
u(z) = h(z), ~ z\in \Gamma,
\end{equation}
where $\Gamma$ is a Jordan curve bounding a Jordan domain $\Omega$ and
$h$ is a prescribed real function.
In this numerical method, first
one computes a AAA approximation $r\approx h$.  This will
normally have poles both inside and outside
$\Gamma$, so its real part, which one might hope would solve the
Laplace problem, will fail to be harmonic in $\Omega$.  The
idea of the AAALS method is to discard the poles of $r$ in $\Omega$ and
use just the real and imaginary parts of the poles outside $\Omega$ to define
a search space for approximating $u$ by a linear least-squares calculation.
The method is surprisingly fast and accurate in many cases, including, for example, on
an L-shaped region as in Figure 4.
So far, there is not much known to explain this success, or to delineate cases
where the AAALS method fails, but
the beginnings of a theory can be found in~\cite{nonconvex}.
The key to taking the theory further probably lies with
Schwarz functions.

The Schwarz function is also relevant to 
analytic continuation of solutions of more general PDE\kern .5pt s problems,
starting from the Hele-Shaw flow problem in fluid
mechanics~\cite{crowdy,gust,howison} and the Laplace equation in 3D~\cite{shapiro}.
Another important case is the
Helmholtz equation $\Delta u + k^2 u = 0$~\cite{barnettb,keller,millar}.
Generalizations of the AAALS method to some of these settings
have proved effective, including
axisymmetric 3D Laplace problems, 
Helmholtz equations~\cite{gt}, and
biharmonic/Stokes flow problems~\cite{xue}.
The present paper is offered partly with these situations
in mind, and also because the topic of Schwarz functions is so interesting
in its own right, so little studied,
and potentially relevant to so many other problems of analytic continuation.

In closing, let us return for a moment to Figure 1.  Comparing
the two sides of the figure, one is reminded that the reflection
across a general arc $\Gamma$  effected by the Schwarz function,
$z \mapsto \overline{S(z)}$, is a generalization of the reflection
in the elementary case in which
$\Gamma$ is a line segment, $z\mapsto \bar z$.  However, in its
impact this is more than ``just'' a generalization because of the
effect visible
repeatedly in our figures: when a boundary $\Gamma$ is curved,
singularities of the Schwarz function
are invariably introduced near $\Gamma$ on the concave side.
These singularities may determine the success or failure of subsequent
numerical operations, and in
particular, as shown in~\cite{nonconvex}, they lead to rational approximations
on nonconvex domains being exponentially more efficient than polynomial ones.
Thus rational approximation is both a good way to calculate Schwarz
functions, and a motivation for wanting to do so.

\section*{Acknowledgments}
I am happy to thank Keaton Burns, Stefano Costa, Kyle
McKee, and Nicholas West for advice on this paper.
Looking back, I am grateful for all I have learned
about rational approximation over the years from Ed Saff.

\end{document}